\documentclass[lettersize,twoside]{IEEEtran}
\usepackage{amsmath,amsfonts}
\usepackage{mathtools}
\usepackage{algorithmic}
\usepackage{algorithm}
\usepackage{array}
\usepackage[caption=false,font=normalsize,labelfont=sf,textfont=sf]{subfig}
\usepackage{textcomp}
\usepackage{stfloats}
\usepackage{url}
\usepackage{svg}
\usepackage{svg-extract}
\usepackage{verbatim}
\usepackage{graphicx}
\usepackage{cite}
\usepackage{hyperref}
\usepackage{tikz}
\usetikzlibrary{arrows.meta, positioning}
\allowdisplaybreaks[3]

\newcommand{\cmmnt}[1]{}
\newcommand{\T}{\mathcal{T}}
\renewcommand{\H}{\mathcal{H}}
\newcommand{\X}{\mathcal{X}}
\begin{document}

\title{A Distributed Optimization Framework to Regulate the Electricity Consumption of a Residential Neighborhood with Renewables~\thanks{Research was partially supported by the DOE under grants DE-EE0009696 and DE-AR-0001282, by the 
NSF under grants CCF-2200052, DMS-1664644, and IIS-1914792, and by the ONR under grant N00014-19-1-2571.}~\thanks{Preprint.}%~\thanks{This work has been submitted to the IEEE for possible publication. Copyright may be transferred without notice, after which this version may no longer be accessible.}
}

%\author{Erhan Can Ozcan, Ioannis Ch. Paschalidis~\IEEEmembership{Fellow,~IEEE,}
\author{%
  \IEEEauthorblockN{%
    Erhan Can Ozcan\IEEEauthorrefmark{2},
    Emiliano Dall'Anese \IEEEauthorrefmark{3} and
    Ioannis Ch. Paschalidis\IEEEauthorrefmark{4},~\IEEEmembership{Fellow,~IEEE}
  }%
  %\IEEEauthorblockA{\IEEEauthorrefmark{1} Affiliation 1}%
  %\IEEEauthorblockA{\IEEEauthorrefmark{2} Affiliation 2}%

        % <-this % stops a space
%\thanks{This paper was produced by the IEEE Publication Technology Group. They are in Piscataway, NJ.}% <-this % stops a space
%\thanks{Manuscript received April 19, 2021; revised August 16, 2021.}}
%\thanks{}}
%\thanks{I. C. Paschalidis is with the Department of Electrical and Computer Engineering, and the Division of System Engineering, Boston University, Boston, MA 02215 USA. e-mail: (yannisp@bu.edu).}
%\thanks{\IEEEauthorrefmark{2} Division of Systems Engineering, Boston University, Boston, MA 02215, e-mail: cozcan@bu.edu.}
\thanks{\IEEEauthorrefmark{2} Erhan Can Ozcan is with Division of Systems Engineering, Boston University, Boston, MA 02215, e-mail: cozcan@bu.edu.}
\thanks{\IEEEauthorrefmark{3} Emiliano Dall'Anese is with Department of Electrical and Computer Eng. and Division of Systems Engineering, Boston University, Boston, MA 02215, e-mail: edallane@bu.edu.}
\thanks{\IEEEauthorrefmark{4} Department of Electrical and Computer Eng., Division of Systems Eng., Dept. of Biomedical Eng., and Faculty of Computing \& Data Sciences, Boston University, 8 St. Mary’s St., Boston, MA 02215, e-mail: yannisp@bu.edu.}}

% The paper headers
%\markboth{IEEE Transactions on Smart Grid, ~Vol.~XX, No.~Y, April~2023}%
%\markboth{}%
%{Ozcan, Dall'Anese and Paschalidis: Distributed Control of Residential Electricity Consumption}

%\IEEEpubid{0000--0000/00\$00.00~\copyright~2021 IEEE}
%\IEEEpubidadjcol
% Remember, if you use this you must call \IEEEpubidadjcol in the second
% column for its text to clear the IEEEpubid mark.

\maketitle
%\IEEEpeerreviewmaketitle

\begin{abstract}

\textcolor{black}{Demand response services at the distribution level are emerging as enabling strategies for improving grid reliability in the presence of intermittent renewable generation and grid congestion.} For residential loads, space heating and cooling, water heating, electric
vehicle charging, and routine appliances make up the bulk of the electricity consumption. \textcolor{black}{Controlling these loads is essential to effectively partake into grid operations and provide services such as peak shaving and demand response.} However, maintaining user comfort is important for ensuring user participation to such a program. This paper formulates a novel mixed integer linear programming problem to control the overall electricity consumption of a residential neighborhood by considering the users' comfort and \textcolor{black}{preferences}. To efficiently solve the problem for communities involving a large number of homes, a distributed optimization framework based on the Dantzig-Wolfe decomposition technique is developed. We demonstrate the load shaping capacity and the computational performance of the proposed optimization framework in a simulated environment.
% Maintaining a balance between electricity demand and supply is a major concern for utilities to ensure a reliable service in power grids, and this problem becomes more acute as the percentage of renewable resources in generation capacity increases due to their intermittent nature. For residential loads, space heating and cooling, water heating, electric
% vehicle charging, and routine appliances make up the bulk of the electricity
% consumption. Controlling these loads can reduce the peak load of a residential
% neighborhood and facilitate matching supply with demand. However, maintaining user
% comfort is important for ensuring user participation to such a program. This paper
% formulates a novel mixed integer linear programming problem to control the overall
% electricity consumption of a residential neighborhood by considering the users'
% comfort. To efficiently solve the problem for communities involving a large number of
% homes, a distributed optimization framework based on the Dantzig-Wolfe decomposition
% technique is developed. We demonstrate the load shaping capacity and the
% computational performance of the proposed optimization framework in a simulated
% environment.

\end{abstract}


\section{Introduction}

% \IEEEPARstart{M}{aintaining} a balance between electricity demand and supply is a major concern for utilities to ensure a reliable service in power grids, especially during peak consumption hours. Although the peak demand occurs occasionally, a significant percentage of the generation capacity is held in reserve to meet it, resulting in increased energy waste \cite{farhangi2009path}. Residential loads, which account for nearly 21\% of the total power consumed in the U.S. \cite{nrel2018}, significantly contribute to daily consumption peaks because residential consumption tends to synchronize with the residents' work schedules ~\cite{he2012residential}.

% \IEEEPARstart{R}{esidential} loads, which account for nearly 21\% of the total power consumed in the U.S. \cite{nrel2018}, significantly contribute to daily consumption peaks because residential consumption tends to synchronize with the residents' work schedules ~\cite{he2012residential}. Although the peak demand occurs occasionally, a significant percentage of the generation capacity is held in reserve to meet it, resulting in increased energy waste \cite{farhangi2009path}.

%OLD VERSION
%%%%%%%
%%%%%%%
%%%%%%%
\IEEEPARstart{E}{nergy} demand in the U.S. grows by almost 1\% annually and is expected to reach approximately 110 quadrillion BTU by 2050~\cite{nalley2022annual}. Renewable energy resources, such as wind turbines and photovoltaic systems, offer significant advantages for society and the environment, and increased production from these sources could accommodate the growing demand. However, the intermittency of renewable resources makes it harder to match supply and demand \cite{luo2015overview}\cmmnt{\cite{uddin2018review, luo2015overview}}. This poses a major concern for utilities striving to ensure a reliable service in power grids, and the vast integration of these resources into power grids requires careful planning \cite{liang2016emerging}\cmmnt{\cite{liang2016emerging, sepasi2023power}}. %This poses a major concern for utilities striving to ensure a reliable service in power grids, especially during peak consumption hours \cite{haider2016review}. Although the peak demand occurs occasionally, a significant percentage of the generation capacity is held in reserve to meet it, resulting in increased energy waste \cite{farhangi2009path}. Residential loads, which account for nearly 21\% of the total power consumed in the U.S. \cite{nrel2018}, significantly contribute to daily consumption peaks because residential consumption tends to synchronize with the residents' work schedules ~\cite{he2012residential}. Furthermore, it is estimated that the energy consumption in residential dwellings will increase by 22\% through 2050, mainly due to population growth ~\cite{nalley2022annual}. Since an increase in the number of end-users can exacerbate the peak load problem in the long run, maintaining a balance between electricity demand and supply can become a more acute problem.

\IEEEpubidadjcol \textcolor{black}{Maintaining a balance between electricity demand and supply becomes particularly critical during peak consumption hours as any imbalance can lead to the failure of the whole grid\cmmnt{\cite{farhangi2009path}}.} \cmmnt{Although the peak demand occurs occasionally, a significant percentage of the generation capacity is held in reserve to meet it, resulting in increased energy waste \cite{farhangi2009path}.} Residential loads, which account for nearly 21\% of the total power consumed in the U.S. \cite{nrel2018}, significantly contribute to daily consumption peaks because residential consumption tends to synchronize with the residents' work schedules\cmmnt{~\cite{he2012residential}}. Furthermore, it is estimated that the energy consumption in residential dwellings will increase by 22\% through 2050, mainly due to population growth ~\cite{nalley2022annual}. Since an increase in the number of end-users can exacerbate the peak load problem in the long run, maintaining a balance between electricity demand and supply can become a more acute problem. A traditional solution to address the peak load problem is to match supply with the load demand by dispatching additional power plants during peak hours \cite{haider2016review}. However, these power plants have high operating costs and can increase pollution, as they primarily rely on fossil fuels \cite{kostkova2013introduction}. Therefore, mechanisms to regulate residential consumption based on available supply and reduce peaks become highly important.
%%%%%%%
%%%%%%%
%%%%%%%

Inducing a change in consumers' energy usage patterns is referred to as {\textit{Demand Response}} (DR)\cmmnt{\cite{albadi2007demand}}. Recent advancements in communication technologies and the proliferation of smart devices have enabled designing {\textit{Home Energy Management Systems}} (HEMS) that can adjust the load consumption schedules of individual homes by considering various criteria, such as \cmmnt{maintaining user comfort, lowering electricity bills, amount of power provided  by renewables, and reducing strain on the power grid by lowering the peak load.}user comfort, electricity rates, the amount of power provided  by renewable resources, and the total strain on the power grid. Therefore, DR strategies have a large potential to make a societal impact.

In this study, we consider a residential neighborhood whose power demand is met by both renewable resources and external traditional generators, and design a DR strategy to match supply with demand while preserving the user comfort. While formulating our problem, we consider various home appliances requiring significant amount of power to function, and model them with enough details to better reflect consumer preferences to make our approach more realistic. Our experiments show that the proposed formulation is capable of controlling the load consumption of the
community based on available supply, and it can mitigate the peak load problem effectively. However, we stress that the
proposed formulation has two major drawbacks that can prevent it from being used in
practice. First, the optimization takes significant amount of time when the number of
participating homes increases since the proposed formulation is a {\em Mixed Integer Linear Programming} (MILP) problem. Even
though there are good commercial solvers to solve MILP problems, they are likely to
fail to provide an effective solution in a reasonable amount of time when the target
community is large. Second, the proposed formulation requires participating homes to
share their preferences and sensitive personal data with the aggregator at
each time interval, which can raise some data privacy concerns among users. In order
to address the aforementioned problems related to efficiency and data privacy, we
develop a distributed optimization framework based on the restricted master heuristic
approach introduced in \cite{my_ref14}. To establish the benefits of the proposed
solution strategy over a commercial solver, Gurobi, we analyze
the solution quality and the optimization time in a simulated
environment. According to our experiments, as the community
size grows, the commercial solver can be ineffective to solve
the centralized problem. However, the proposed distributed optimization framework can provide an effective solution with a small optimality gap in significantly less amount
of time even if the number of participating homes becomes
large (e.g., 10, 000).

%Load shifting, which reschedules some appliances to run during off-peak hours, is one of the oldest DR strategies to reduce the peak load. However, it may result in the relocation of the existing peak load unless the load is distributed over several off-peak time periods \cite{liefficient2017}. Furthermore, 

\subsection{Related Work}

Continuous user participation is essential for the success of residential DR programs; thus, maintaining the comfort of participants is vital while adjusting their load schedules \cite{shimomura2014method}. There are numerous studies where the user comfort is represented by utility functions \cite{yu2015real, jiang2020residential, wen2022demand}. However, as there is no simple way to translate user preferences into utility functions, it becomes hard to assess whether the scheduled load plan meets the user preference or not. On the contrary, {\textit{Model Predictive Control}} (MPC)-based approaches can ensure user comfort by incorporating inputs such as user preferences, weather forecasts, appliance properties, and appliance dynamics into the decision-making stage \cite{parisio2015mpc, foresee2017}, thus facilitating the design of effective DR programs for the residential sector. Therefore, we design an MPC-based DR strategy considering various home appliances in our paper. \cmmnt{Considering the success of MPC-based DR strategies in some recent studies ~\cite{bandyopadhyay2020one, chen2023multi}, we design an MPC-based DR strategy by considering various home appliances in this paper.}

Although various MPC-based DR programs have been proposed in the literature, some earlier works, such as \cite{du2011appliance,chen2013mpc, foresee2017, li2020real}\cmmnt{\cite{du2011appliance,chen2013mpc, foresee2017, bandyopadhyay2020one, li2020real}} focus on individual homes without formulating an optimization problem at the neighborhood level. These strategies cannot ensure coordination toward a specific goal, such as evenly distributing the power across the planning horizon or adapting consumption based on the available supply in the grid. \cmmnt{Controlling the load consumption of a community requires aggregating many participants, and coordinating their actions based on their preferences. In this regard, \cite{chapman2016algorithmic} provides a framework for designing a DR strategy that facilitates coordination among participants, and presents the practical challenges that may arise due to different design choices. However, we will explicitly review the existing studies  some common}Controlling the load consumption of a community requires aggregating many participants, and a general framework is provided in \cite{chapman2016algorithmic} for designing a DR strategy that facilitates coordination among participants. In the literature, a large group of studies focus on controlling the peak demand. For example, \cite{zhu2012integer} formulates\cmmnt{studies in \cite{zhu2012integer, shakouri2017multi} formulate} an optimization problem at the neighborhood level to minimize the peak load consumption, while studies in \cite{meng2013stackelberg, li2021coordinating, rezaei2022hierarchical, mhanna2016fast}\cmmnt{\cite{meng2013stackelberg, carrasqueira2017bi, li2021coordinating, rezaei2022hierarchical, mhanna2016fast}} restrict the peak load consumption by adding a hard constraint on the overall load consumption of the community. While these strategies can be effective in controlling peak demand, fluctuations in load consumption may still occur over the planning horizon. This can pose a significant challenge in energy generation planning, especially while managing a large percentage of renewable resources in the generation mix, as it is necessary to maintain a balance between supply and demand for reliable grid service.

On the other hand, studies in \cite{safdarian2015optimal} and \cite{chen2023multi} propose DR strategies that keep the aggregated load consumption of the community close to a target level rather than minimizing the peak load. Despite testing their models under various pricing structures, these two studies simply use existing market prices without optimizing them. However, each home has different consumption habits and each home may respond differently to price changes. Therefore, these approaches may not realize their full potential. The authors in \cite{ozcan2023stackelberg} propose a Stackelberg game between a coordination agent and participating homes. While each home optimizes a comfort-cost trade off to determine a load schedule of its available appliances in response to a price vector set by the coordination agent, the goal of the coordination agent is to find a price vector that will keep the overall load consumption of the neighborhood around a predetermined target value. However, the resulting optimization problem is non-convex and the proposed gradient-based solution approach can be susceptible to local minima. \cmmnt{The studies presented in \cite{safdarian2015optimal, chen2023multi, ozcan2023stackelberg} assume that the main goal of the participating homes is to minimize the electricity bill, thus designing DR strategies based on this assumption.}The studies presented in \cite{safdarian2015optimal, chen2023multi, ozcan2023stackelberg} designed residential DR strategies assuming that the main goal of the participating homes is to minimize the electricity bill. However, unlike the industrial and commercial sectors, where profit maximization is a driving factor, maintaining the comfort level can be more important than lowering the electricity bill in the residential sector \cite{yan2018review}. Furthermore, \cite{reiss2005household} show that customers might be non-responsive to the electricity price. \textcolor{black}{In this regard, we design a DR strategy, in which incentive payments are distributed to the participating homes to keep the aggregated load consumption of the community close to a target level.} 

Several studies have proposed residential DR strategies by integrating renewable energy resources. For example, the studies in \cite{liu2018day,kim2013scalable,xu2016hierarchical,park2017residential} formulate a problem to minimize the cost of purchasing electricity from external power systems without focusing on how the purchased power changes throughout the planning horizon. \cmmnt{\textcolor{red}{In \cite{xu2016hierarchical}, a two-stage stochastic problem is proposed to handle uncertainties in the electricity price and renewable generation. However, similar to \cite{liu2018day,kim2013scalable} and \cite{park2017residential}, the primary goal of this strategy is the cost minimization without }proposes a two-stage stochastic problem to handle uncertainties in the electricity price and renewable generation for cost minimization.} In \cite{zhang2019distributionally}, a robust optimization problem is designed to utilize all the power generated by renewable resources while keeping the room temperature within desirable limits. As long as renewable resources provide sufficient power to meet the overall demand, this formulation is applicable. However, it remains unclear how to effectively manage load consumption in case of excessive demand. In \cite{khezri2021demand}, a demand side management strategy is proposed to promote lowering electricity consumption when generation from renewable resources is scarce. However, it does not take into account how user comfort is affected by the delay or curtailment of loads, thus poorly representing the user aspect. 

\textcolor{black}{Developing practical and efficient solution mechanisms is as important as designing effective DR strategies. This study formulates an MILP problem at the neighborhood level to control the load demand, and MILPs are non-convex and known to be NP-hard\cite{testa2019distributed, camisa2021distributed}. Even though commercial solvers efficiently handle small-scale MILP problems using strategies based on branch-and-bound \cite{achterberg2013mixed,morrison2016branch}, problems related to residential DR are generally in large-scale. As a result, the performance of off-the-shelf solvers can be limited, making it necessary to develop efficient solution strategies that are feasible to implement in practice. For example, \cite{kim2013scalable} proposes a distributed solution strategy based on Lagrangian relaxation to solve a demand response problem with mixed-integer constraints. \textcolor{black}{However, this type of strategies can suffer from feasibility issues as indicated in \cite{kim2013scalable}.} Similarly, some recent studies propose utilizing \textit{Alternating Direction Method of Multipliers (ADMM)} to solve MILP problems \cite{my_ref_ADMM_1}, \cite{my_ref_ADMM_2}. However, convergence guarantees to global optimal solutions have not been established for MILP problems using ADMM \textcolor{black}{unless the problem admits a certain block-angular structure} \cite{sun2024decomposition}. Furthermore, the performance of ADMM in solving large scale MILP problems has not been investigated in depth, yet. \textcolor{black}{Finally, due to the objective function selected in our formulation, decomposing the problem via ADMM is only possible after fixing the load consumption of other homes at each iteration, which can slow down the optimization.}} In our paper, we develop a distributed optimization framework based on the restricted master heuristic approach introduced in \cite{my_ref14}, \textcolor{black}{and employed a delayed column generation strategy to ensure feasibility.}

\subsection{Contributions}
Our main contributions are as follows: $(1)$ We design a DR strategy that maintains a balance between supply and demand, thus integrating renewable resources into the generation mix effectively. Based on our experiments, the formulated MILP problem can help reduce peak load while controlling the activity of various appliances based on various comfort-related constraints. $(2)$ We show that the formulated optimization problem can be decomposed across homes, and develop an approach based on the Dantzig-Wolfe decomposition technique to solve the problem in a distributed way. \textcolor{black}{Therefore, we propose a distributed algorithm that can be implemented into existing off-the-shelf solvers in HEMS.} $(3)$ We compare the performance of the proposed approach with a commercial solver trying to solve the centralized version of the problem, and demonstrate that our algorithm can provide an effective solution \textcolor{black}{(with optimality gap less than 1\%)} in a short period of time.

\textcolor{black}{\textbf{Notation:} We denote vectors by bold lowercase letters (e.g., $\boldsymbol{x}$), and represent the $k$-th element of the vector $\boldsymbol{x}$ by $x(k)$. Bold uppercase letters (e.g., $\boldsymbol{A}$) and calligraphic letters (e.g., $\mathcal{P}$) are used to denote matrices and sets, respectively. We represent the comfort of participants by a mixed-integer polyhedron, of which the general form is: \[\mathcal{P}=\{(\boldsymbol{x},\boldsymbol{y})\ |\ \boldsymbol{A}\boldsymbol{x}+\boldsymbol{B}\boldsymbol{y}\leq \boldsymbol{b},\ (\boldsymbol{x},\boldsymbol{y})\in \mathbb{Z}^{d_x}\times \mathbb{R}^{d_y}\},\]
where $\boldsymbol{A}$ is the constraint matrix, $\boldsymbol{b}$ is the right hand side vector, and $(\boldsymbol{x}, \boldsymbol{y})$ is the decision vector with $d_x$ integer components and $d_y$ real components, respectively.}

\section{Proposed Optimization Problem}

Consider a power grid consisting of an aggregator and $N$ participating homes,
each equipped with controllable $M$ appliances \footnote{For simplicity of exposition and to maintain the notation light, we assume that each home has M controllable loads. However, our approach can be of course applied to the case where each home has a different numbers of controllable loads.}. In our study, we assume that load consumption associated with heating, ventilation, and air-conditioning system (HVAC), electric vehicle (EV), electric water heater (EWH), and basic appliances such as washing machine, dryer, and oven can be controlled. Suppose that $p_{ij}(t)$ is a decision variable denoting the load consumption of home $i$ for controllable appliance $j$ at time point $t$, and \[
\boldsymbol{p}_i=(p_{i1}(0) \ldots p_{i1}(K-1) \ldots p_{iM}(0) \ldots p_{iM}(K-1))
  \]
is a vector of decision variables denoting the load schedule of
home i for M controllable loads over the next K time intervals. The comfort related constraints of each home are designed by considering the individual preferences of the corresponding home for these appliances, and these constraints form a mixed-integer polyhedron $\H_i$, of which general form is as follows:
\begin{equation*}
\begin{split}
\mathcal{H}_i=\{(\boldsymbol{p}_i,\boldsymbol{y}_i, \boldsymbol{s}_i)\ |\ &\boldsymbol{A}_i\boldsymbol{p}_i+\boldsymbol{B}_i\boldsymbol{y}_i + \boldsymbol{C}_i\boldsymbol{s}_i\leq \boldsymbol{b}_i,\ \\&(\boldsymbol{p}_i,\boldsymbol{y}_i,\boldsymbol{s}_i)\in \mathbb{R}^{d_p}\times\mathbb{Z}^{d_y}\times \mathbb{R}^{d_s}\},
\end{split}
\end{equation*}
where $\boldsymbol{p}_i$ is the load consumption vector with $d_p$ real components; $\boldsymbol{y}_i$ and $\boldsymbol{s}_i$ are auxiliary variables with with $d_y$ integer components and $d_s$ real components, respectively; $\boldsymbol{A}_i$, $\boldsymbol{B}_i$, and $\boldsymbol{C}_i$ are the constraint matrices defining the comfort of the user; and $\boldsymbol{b}_i$ is the right hand side vector. To keep notation simpler, we omit the auxiliary variables, and represent the user comfort by $\boldsymbol{p}_i\in\mathcal{H}_i$ throughout the text. The explicit definition of the constraints are available in Appendix \ref{app_homeappliance}.\cmmnt{of which explicit definition is moved to Appendix \ref{app_homeappliance} to keep notation simpler. Hence, user comfort can be ensured as long as $\boldsymbol{p}_i\in\H_i, \ \forall i$.} However, each home $i$ may have a desirable load consumption
schedule due to some specific consumption habits, \textcolor{black}{and users can be reluctant to deviate from their desirable load schedules. Therefore, the aggregator distributes incentive payments to the participating homes.} We denote by $\overline{\boldsymbol{p}_i}$ the
desirable load consumption schedule of home $i$ and by $\boldsymbol{u_i}$ the
deviations such that:
\[
p_{ij}(t)=\overline{p_{ij}}(t) + u_{ij}(t)\quad\forall\ i,j,t.
\]

\textcolor{black}{Apart from the controllable loads, each home has some uncontrollable load consumption, which we denote by $p_{i,\text{U}}(t)$ for home $i$; thus, the overall load demand at time $t$ is equal to $\sum_{i=1}^N\sum_{j=1}^M p_{ij}(t)\ +\ \sum_{i=1}^Np_{i,\text{U}}(t)$. To meet the load demand, both renewable resources and external traditional generators are utilized, and the aggregator controls the available supply in the grid by providing additional power from external generators. In other words, we assume that renewable energy resources provides power by $p^{r}(t)$ at time $t$, and the aggregator adjusts the supply using additional power sources to meet the demand.}

\textcolor{black}{In our formulation, the role of the aggregator is to maintain a balance between supply and demand while minimizing the total incentive payments distributed to the community. Therefore, we formulate the following optimization problem:}
%
%To mitigate the peak load problem, we suggest to keep 
%
% To satisfy the power demand during the planning horizon, the coordination agent needs to purchase some power from the electricity market. Let denote the amount of power that CA buys from the power market, and we assume that the goal is to keep the purchased power as flat as possible during the planning horizon.
%
%The role of the aggregator is to provide enough power to the grid
\begin{subequations}
\begin{align}
\min_{\boldsymbol{p}_1,\cdots,\boldsymbol{p}_N,\atop \boldsymbol{u}_1,\cdots,\boldsymbol{u}_N, q, \boldsymbol{a}} &\sum_{t=0}^{K-1} \left|a(t)\right| + \sum_{i=1}^N \sum_{j=1}^M\sum_{t=0}^{K-1}c_{ij}\left|u_{ij}(t)\right| \label{eq_4}  \\
\text{s.t.} \quad &p_{ij}(t)=\overline{p_{ij}}(t) +u_{ij}(t) \label{eq_dev},\quad \quad
\ \ \forall\ i,j,t, \\
&a(t)=q + p^{r}(t) - \sum_{i=1}^N\sum_{j=1}^M p_{ij}(t) \notag \\
&\  \qquad \qquad \qquad\qquad  -\sum_{i=1}^Np_{i,\text{U}}(t),\quad \forall\ t, \label{eq_mismatch} \\
%&q(0) = q(1) = \cdots = q(K-1)\label{eq_qbalance}, \\
%&q(t) \geq 0, \qquad\qquad\qquad\qquad\quad \forall\ t \label{eq_qnonnegative}, \\
&q \geq 0, \qquad\qquad\qquad\qquad\quad \label{eq_qnonnegative} \\
&\boldsymbol{p}_i \in \H_i,\quad\quad\quad\quad\quad\quad\quad\quad\quad\ \ \forall\ i, \label{eq_comfort} 
\end{align}
\end{subequations}
% \textcolor{red}{where $\boldsymbol{q}$ is the decision vector denoting the amount of power that the aggregator provides from external generators during the planning horizon, $c_{ij}$ represents the rate of incentive payment paid for appliance $j$ for home $i$, and $|\cdot|$ denotes the absolute value function.} The constraints in Equations \eqref{eq_dev} and \eqref{eq_mismatch} measure the deviation from desirable load level and the mismatch between supply and demand, respectively. Equation \eqref{eq_qbalance} ensures that the power provided by external generators remains constant throughout the planning horizon to prevent instantaneous peak load consumption. Finally, the constraint in \eqref{eq_comfort} maintains the user comfort. On the other hand, the first term in the objective function helps the aggregator to minimize the mismatch between supply and demand to keep the load consumption of the community under control, while the second term minimizes the incentive payments distributed to the community. Lastly, since the comfort related constraints of each home are represented by a mixed integer polyhedron, the optimization problem becomes MILP.
\textcolor{black}{where $q$ denotes the amount of power that the aggregator provides from external generators at each time period $t$}, $c_{ij}$ represents the rate of incentive payment paid for appliance $j$ for home $i$, and $|\cdot|$ denotes the absolute value function. The constraints in Equations \eqref{eq_dev} and \eqref{eq_mismatch} measure the deviation from desirable load level and the mismatch between supply and demand, respectively. \textcolor{black}{In our formulation, the power provided by external generators, $q$, remains constant throughout the planning horizon to prevent instantaneous peak load consumption.} Finally, the constraint in \eqref{eq_comfort} maintains the user comfort. On the other hand, the first term in the objective function helps the aggregator to minimize the mismatch between supply and demand to keep the load consumption of the community under control, while the second term minimizes the incentive payments distributed to the community. Lastly, since the comfort related
constraints of each home are represented by a mixed integer polyhedron, the
optimization problem becomes MILP.

\section{Optimization Framework}
\label{section:opt}

\textcolor{black}{Although the centralized problem can be optimized by some commercial solvers as it
stands, it is possible to decompose the problem at the home level after some
modifications. Hence, we replace the absolute values in the objective function by an
equivalent reformulation, and we define a new set $\X_i$ by combining the constraints \eqref{eq_dev} and \eqref{eq_comfort}. The resulting optimization problem is as
follows:}

\noindent\textbf{Centralized-IP:}
\begin{subequations}
\label{central_ip}
\begin{align}
\min_{\boldsymbol{p}_1,\cdots,\boldsymbol{p}_N,\atop \boldsymbol{u}^+_1,\cdots,\boldsymbol{u}^+_N, q,\boldsymbol{s},\boldsymbol{a}} &\sum_{t=0}^{K-1} s(t) + \sum_{i=1}^N \sum_{j=1}^M\sum_{t=0}^{K-1}c_{ij} u_{ij}^+(t) \label{eq_8}  \\
\text{s.t.} \quad &s(t) - a(t) \geq 0,\quad\quad\quad\quad\quad\quad\ \forall\ t,\\
&s(t) + a(t) \geq 0,\quad\quad\quad\quad\quad\quad\ \forall\ t,\\
&a(t) - q +\sum_{i=1}^N\sum_{j=1}^M p_{ij}(t) \notag \\
&\quad\quad = p^{r}(t) - \sum_{i=1}^Np_{i,\text{U}}(t),\quad\quad \forall t,\\
%&q(0) = q(1) = \cdots = q(K-1), \\
&q \geq 0, \qquad\qquad\qquad\qquad\quad \\
&\left(\boldsymbol{p}_i,\boldsymbol{u}^+_i\right) \in
\X_i,\quad\quad\quad\quad\quad\quad\ \ \ \forall\ i,  \label{eq_cent_last}
\end{align}
\end{subequations}
% \noindent\textbf{Centralized-IP:}
% \begin{align}
% \min_{\boldsymbol{p}_1,\cdots,\boldsymbol{p}_N,\atop \boldsymbol{u}^+_1,\cdots,\boldsymbol{u}^+_N} &\sum_{t=0}^{K-1} s(t) + \sum_{i=1}^N \sum_{j=1}^M\sum_{t=0}^{K-1}c_{ij} u_{ij}^+(t) \label{eq_8}  \\
% \text{s.t.} \quad &s(t) - a(t) \geq 0,\quad\quad\quad\quad\quad\quad\ \forall\ t,\notag\\
% &s(t) + a(t) \geq 0,\quad\quad\quad\quad\quad\quad\ \forall\ t,\notag\\
% &a(t)+\sum_{i=1}^N\sum_{j=1}^M p_{ij}(t)=Q(t),\quad \forall\ t,\notag\\
% &\left(\boldsymbol{p}_i,\boldsymbol{u}^+_i\right) \in
% \X_i,\quad\quad\quad\quad\quad\quad\ \ \ \forall\ i, \notag
% \end{align}
where
\begin{equation*}
\begin{split}
\X_i = \{\left(\boldsymbol{p}_i,\boldsymbol{u}^+_i\right)\  |\ &\boldsymbol{p}_i \in
\H_i; \notag \\
& p_{ij}(t)=\overline{p_{ij}}(t) +u_{ij}(t),\  \forall\  j,t;\notag\\
& u^+_{ij}(t) = |u_{ij}(t)|,\  \forall\ j,t\}. \notag
\end{split}
\end{equation*}
In order to solve the optimization problem given in \eqref{eq_8}\textendash\eqref{eq_cent_last}, we employ a
distributed optimization framework. 
%based on the restricted master heuristic approach introduced
%in \cite{my_ref14}. Therefore, we use a similar terminology while discussing the
%details of our approach.

The \textcolor{black}{set} $\X_i$ is bounded for all homes, and suppose that $G_i$ is the
generating set that consists of the feasible solutions to $\X_i$, i.e.,
\[
\X_i=\left\{\boldsymbol{x}_i^g:\ \boldsymbol{x}_i^g=\left(\boldsymbol{p}_i^g,
\boldsymbol{u}_i^{+,g}\right) \right\}_{g \in G_i}.
\]
Although a {\em Master Problem (MP)} equivalent to Centralized-IP \textcolor{black}{in \eqref{central_ip}} can be formulated
by using the generating set $G_i$, optimizing the formulated MP can be time consuming
as the the number of feasible points grows exponentially \textcolor{black}{\cite{puchinger2011dantzig,el2023consensus}}. Hence, we define a {\em
  Restricted Master Problem (RMP)} by considering only a small subset of feasible
points in $\X_i$, and gradually increase the problem complexity by using a column
generation idea.

Suppose that $\widehat{G_i}$ includes only a small subset of feasible points in $\X_i$, e.g., $\widehat{G_i} \subset G_i$, then the RMP can be formulated as follows:
\begin{subequations}
\label{prob_rmp}
\begin{align}
\notag \textbf{RMP:}\\
\min_{\lambda_i^g, q,\boldsymbol{s},\boldsymbol{a}} &\sum_{t=0}^{K-1} s(t) + \sum_{i=1}^N\sum_{g \in \widehat{G_i}}\left(\sum_{j=1}^M\sum_{t=0}^{K-1} c_{ij} u_{ij}^{+,g}(t)\right)\lambda_i^g\label{eq_13}  \\
\text{s.t.} \quad &s(t) - a(t) \geq \label{eq_rmp_first}
0,\quad\quad\quad\quad\quad\quad\quad\quad\quad\quad\quad \forall\ t, \\
&s(t) + a(t) \geq 0,\quad\quad\quad\quad\quad\quad\quad\quad\quad\quad\quad
\forall\ t, \\
&a(t)-q \notag\\
&\qquad+\sum_{i=1}^N\sum_{g \in \widehat{G_i}}\left(\sum_{j=1}^M
p_{ij}^g(t)\right)\lambda_i^g\notag\\
&\qquad=p^{r}(t) - \sum_{i=1}^Np_{i,\text{U}}(t),\ \qquad\qquad \forall\ t, \\
&\sum_{g \in
  \widehat{G_i}}\lambda_i^g=1,\quad\quad\quad\quad\quad\quad\qquad\quad\quad \forall\ i,\label{rmp_eq_lamdasum}\\
%&q(0) = q(1) = \cdots = q(K-1), \label{eq_rmp_last} \\
&q \geq 0, \qquad\qquad\qquad\qquad\quad\qquad\quad \\
&\lambda_i^g \in \{0,1\},\qquad\qquad\qquad\qquad\quad \forall i, \forall g\in \widehat{G_i},\label{rmp_last}
\end{align}
\end{subequations}
% \begin{align}
% \notag \textbf{RMP:}\\
% \min_{\lambda_i^g} &\sum_{t=0}^{K-1} s(t) + \sum_{i=1}^N\sum_{g \in \widehat{G_i}}\left(\sum_{j=1}^M\sum_{t=0}^{K-1} c_{ij} u_{ij}^{+,g}(t)\right)\lambda_i^g\label{eq_13}  \\
% \text{s.t.} \quad &s(t) - a(t) \geq
% 0,\quad\quad\quad\quad\quad\quad\quad\quad\quad\quad\quad \forall\ t,\notag \\
% &s(t) + a(t) \geq 0,\quad\quad\quad\quad\quad\quad\quad\quad\quad\quad\quad
% \forall\ t, \notag\\
% &a(t)+\sum_{i=1}^N\sum_{g \in \widehat{G_i}}\left(\sum_{j=1}^M
% p_{ij}^g(t)\right)\lambda_i^g=Q(t),\ \ \forall\ t, \notag \\
% &\sum_{g \in
%   \widehat{G_i}}\lambda_i^g=1,\quad\quad\quad\quad\quad\quad\qquad\quad\quad\quad\quad\ \forall\ i,
% \notag\\
% &\lambda_i^g \in \{0,1\},\qquad\qquad\qquad\qquad\quad \forall i, \forall g\in \widehat{G_i},\notag
% \end{align}
where $\lambda_i^g$ is the variable weight of the associated feasible point,
$\boldsymbol{x}_i^g$, of home $i$. By removing the integrality constraints on
$\lambda_i^g$ variables in the RMP, we obtain the relaxed-RMP, and we use the dual
information provided by the relaxed-RMP in our column generation mechanism. The
relaxed-RMP is available in Appendix \ref{app_distopt}.

Suppose that the relaxed-RMP has dual solution vectors
$\boldsymbol{\sigma}_1,\boldsymbol{\sigma}_2,\boldsymbol{\sigma}_3,$, and $\boldsymbol{\sigma}_4$ associated with the constraints in \eqref{eq_rmp_first}\textendash\eqref{rmp_eq_lamdasum} in the
order they appear. Then, each home $i$ can choose which feasible point to share with
the aggregator by optimizing the following pricing subproblem:
\begin{equation}
\begin{split}
&L_i(\boldsymbol{\sigma}_1,\boldsymbol{\sigma}_2,\boldsymbol{\sigma}_3,\boldsymbol{\sigma}_4):= \\
&\min_{(\boldsymbol{p}_i, \boldsymbol{u}^+_i) \in
    X_i}\left(\sum_{j=1}^M\sum_{t=0}^{K-1} c_{ij}
  u_{ij}^{+}(t)-\sum_{t=0}^{K-1}\sum_{j=1}^M p_{ij}(t)\sigma_3(t)\right). \label{eq_19}
\end{split}
\end{equation}

\textcolor{black}{Similar to Problem in \eqref{central_ip}}, the formulated pricing problem is optimized subject to the
constraints in $\X_i$; thus, it is a MILP problem. However, since the pricing problem
has only a few binary variables at the home level, any commercial solver embedded in
HEMS can optimize it efficiently, and HEMS decides whether a new feasible point is
shared with the aggregator according to the following simple rule.  If
$L_i(\boldsymbol{\sigma}_1, \boldsymbol{\sigma}_2, \boldsymbol{\sigma}_3,
\boldsymbol{\sigma}_4)<\boldsymbol{\sigma}_4(i)$, $\boldsymbol{x}_i^*$, the optimal
solution of the optimization problem given in (\ref{eq_19}) is sent to the
aggregator, and it is added to the subset $\widehat{G_i}$ since the reduced
cost of this new feasible point is negative in the relaxed-RMP. Otherwise,
$\widehat{G_i}$ remains unchanged in the corresponding column generation iteration.

Column generation iterations can continue until
$L_i(\boldsymbol{\sigma}_1,\boldsymbol{\sigma}_2,\boldsymbol{\sigma}_3,\boldsymbol{\sigma}_4)\geq\boldsymbol{\sigma}_4(i),
\ \forall i.$ However, in practice, it might take substantial amount of time to
satisfy this stopping condition \cite{my_ref15}. We know that the
improvement in the objective value of the relaxed-RMP slows down gradually, and we
can cease the column generation iterations when the current cost, $z^{rRMP}$, is
close enough to the Lagrangian dual bound \cite{my_ref15}. Fortunately, the
Lagrangian dual bound is readily available once the subproblems are
optimized\cite{my_ref16}, and it can be calculated as follows:
\begin{equation}
\begin{split}
\xi
(\boldsymbol{\sigma}_1,\boldsymbol{\sigma}_2,\boldsymbol{\sigma}_3,\boldsymbol{\sigma}_4)&=\sum_{t=0}^{K-1}(p^{r}(t) - \sum_{i=1}^Np_{i,\text{U}}(t))\sigma_3(t) \\ &+\sum_{i=1}^NL_i(\boldsymbol{\sigma}_1,\boldsymbol{\sigma}_2,\boldsymbol{\sigma}_3,\boldsymbol{\sigma}_4).  \label{eq_20}
\end{split}
\end{equation}
In \cite{my_ref17}, it is discussed that optimizing the pricing problem in a more
restricted region can yield stronger dual bounds at each iteration, thereby
terminating the CG iterations earlier. The illustrative example provided in Appendix
\ref{app_distopt} shows that our approach \textcolor{black}{will} indeed enjoy stronger dual bounds.

Even though we use the revised simplex to optimize the relaxed-RMP after adding new
columns, optimization time increases gradually as the model complexity grows. Hence,
we remove some of the feasible points that are added to $\widehat{G_i}$ if their
corresponding coefficients remain zero for $\kappa$ consecutive iterations. Although
removing the unused feasible points increases the number of CG iterations that our
algorithm needs before termination, it also decreases the optimization time of the
relaxed-RMP considerably, which helps us to reduce the total optimization time.

Let $\lambda_i^{g,rRMP}$ be the optimal solution of the relaxed-RMP for home $i$,
$\forall\ g \in \widehat{G_i}$ after the column generation iterations are
completed. Then, the resulting optimal solution $\boldsymbol{x}_{i}^{rRMP}$ can be
expressed as follows:
\[
\boldsymbol{x}_{i}^{rRMP}=\sum_{g \in
  \widehat{G_i}}\boldsymbol{x}_i^g\lambda_i^{g,rRMP},\ \forall\ i.
\]

Unfortunately, the resulting optimal solution of some homes in the community can be
infeasible, e.g., $\boldsymbol{x}_{i}^{rRMP} \notin \X_i,$ for some $i$'s, and we
need to take an additional step to recover feasible solutions after solving the
relaxed-RMP. Since we solve the pricing problem by considering the subsystem $\X_i$,
each point in the subset $\widehat{G_i}$ is a valid solution for home $i$. Therefore,
we reformulate the RMP given in \eqref{prob_rmp} by considering the learned subsets
$\widehat{G_i},\ \forall\ i$; and optimize as the last step.

The pseudocode showing the aforementioned steps of the distributed optimization
algorithm is given in Algorithm \ref{alg:one}. While the outer loop in Algorithm
\ref{alg:one} controls the column generation iterations by checking the optimality gap
at each iteration, the inner loop finds alternative consumption schedules for each
home. At this point, we want to emphasize that each home performs the inner loop
steps independently, and thus, the designed algorithm solves the issues related to
efficiency and data privacy. \textcolor{black}{Figure \ref{fig:alg} shows the information exchange between participating homes and the aggregator.}

\begin{figure}
    \centering
    \begin{tikzpicture}[
  node distance=1.8cm and 3cm,
  every node/.style={align=center},
  box/.style={rectangle, draw=black, minimum width=2.0cm, minimum height=0.8cm},
  widebox/.style={rectangle, draw=black, minimum width=7.5cm, minimum height=1.2cm},
  dbox/.style={rectangle, draw=black, dashed, minimum width=6cm, minimum height=1.2cm},
  arrow/.style={-{Triangle}, thin, black},
  doublearrow/.style={<->, thick, blue},
  dashedarrow/.style={-{Triangle}, thin, black, dashed}
  %arrowbox/.style={rectangle, draw=black, width=1.0cm height=1.0cm}
]

% Top dashed box
\node[dbox] (structured) at (0,5.5) {
1. Solve $L_i(\boldsymbol{\sigma}_1,\boldsymbol{\sigma}_2,\boldsymbol{\sigma}_3,\boldsymbol{\sigma}_4)$ to obtain $\boldsymbol{x}_i$.\ \ \ \\ 
2. If $L_i(\boldsymbol{\sigma}_1,\boldsymbol{\sigma}_2,\boldsymbol{\sigma}_3,\boldsymbol{\sigma}_4) \geq \boldsymbol{\sigma}_4(i)$, send $\boldsymbol{x}_i$.};

% Middle row boxes
\node[box] (social) at (-3.,3.75) {Home 1};
\node[box] (crm) at (0,3.75) {Home $i$};
\node[box] (erp) at (3,3.75) {Home $N$};

% Data Lake center
\node[widebox] (lake) at (0,1.75) {Aggregator};

% Data Science Workbench
\node[dbox, minimum width=5.5cm] (workbench) at (0,-0.) {Solve relaxed-RMP to obtain $\boldsymbol{\sigma}$.};

% Output box
%\node[dbox, minimum width=10cm] (output) at (0,-3.5) {Dashboards / Reports / Machine Learning Models};

% Arrows from structured to CRM
\draw[dashedarrow] (crm) -- (structured);

% Bidirectional arrows from sources to lake
\draw[arrow] ([xshift=-6.pt]social.south) -- ([xshift=15pt]lake.north west);
\draw[arrow] ([xshift=27pt]lake.north west) -- ([xshift=6pt]social.south);
\node at (-3.5, 2.75) {$\boldsymbol{x}_1$};
\node at (-2.5, 2.75) {$\boldsymbol{\sigma}$};

\draw[arrow] ([xshift=-6pt]crm.south) -- ([xshift=-6pt]lake.north);
\draw[arrow] ([xshift=6pt]lake.north) -- ([xshift=6pt]crm.south);
\node at (-0.5, 2.75) {$\boldsymbol{x}_i$};
\node at (0.5, 2.75) {$\boldsymbol{\sigma}$};

%\draw[arrow] ([xshift=-6pt]erp.south) -- ([xshift=-6pt]lake.north east);
%\draw[arrow] ([xshift=6pt]lake.north east) -- ([xshift=6pt]erp.south);
\draw[arrow] ([xshift=-6pt]erp.south) -- ([xshift=-27pt]lake.north east);
\draw[arrow] ([xshift=-15pt]lake.north east) -- ([xshift=6pt]erp.south);
\node at (2.5, 2.75) {$\boldsymbol{x}_N$};
\node at (3.5, 2.75) {$\boldsymbol{\sigma}$};

% Dashed arrow from lake to workbench
\draw[dashedarrow] (lake.south) -- (workbench.north);

% Dashed arrow from workbench to output
%\draw[dashedarrow] (workbench.south) -- (output.north);

\end{tikzpicture}
    \caption{Information exchange between participating homes and the aggregator.}
    \label{fig:alg}
\end{figure}

% \begin{figure}[htb]
% \centering
% \includegraphics[width=2.75in]{pseudo.png}
% %\vspace{-1.5cm}
% \caption{Information exchange between participating homes and the aggregator.} 
% \label{fig:alg}
% \end{figure}

\begin{algorithm}[ht]
\caption{Restricted Master Heuristic with a given $\epsilon>0$, and integer $\kappa$
  as a hyperparameter.} 
\label{alg:one}
\begin{algorithmic}[1]
\STATE \textbf{Initialize}
\STATE Set $\widehat{G_i}=\emptyset,\ \forall\ i.$
\STATE Each home has a feasible vector, $\boldsymbol{x}_i \in \X_i$, ready to send. 
\REPEAT 
\STATE Receive the vector $\boldsymbol{x}_i,\ \forall\ i.$
\STATE $\widehat{G_i}=\widehat{G_i}\cup \boldsymbol{x}_i,\ \forall\ i.$
\STATE Solve the relaxed-RMP to obtain $\boldsymbol{\sigma}_1,\boldsymbol{\sigma}_2,\boldsymbol{\sigma}_3,\boldsymbol{\sigma}_4$, $z^{rRMP}$, and $\lambda_i^{g,rRMP}$.
\STATE i=1
\REPEAT 
\STATE Solve $L_i(\boldsymbol{\sigma}_1,\boldsymbol{\sigma}_2,\boldsymbol{\sigma}_3,\boldsymbol{\sigma}_4)$ to obtain $\boldsymbol{x}_i^*.$
\IF {$L_i(\boldsymbol{\sigma}_1,\boldsymbol{\sigma}_2,\boldsymbol{\sigma}_3,\boldsymbol{\sigma}_4)<\boldsymbol{\sigma}_4(i)$}
\STATE $\boldsymbol{x}_i=\boldsymbol{x}_i^*$ %\hfill\COMMENT{\forall i}
\ELSE
\STATE $\boldsymbol{x}_i=\emptyset$.
\ENDIF
\STATE $G^\lambda=\{g\ |\ \lambda_i^{g,rRMP}=0 \text{ for the last } \kappa\ \text{iterations}\}$
\STATE $\widehat{G_i}=\widehat{G_i}\setminus \boldsymbol{x}_i^g,\ \forall\  g \in G^\lambda$.
\STATE $i=i+1$
\UNTIL {$i>N$.}
\STATE Calculate $\xi
(\boldsymbol{\sigma}_1,\boldsymbol{\sigma}_2,\boldsymbol{\sigma}_3,\boldsymbol{\sigma}_4)$
as in (\ref{eq_20}).
\UNTIL{$\ \ \frac{\left|z^{rRMP}\ - \ \xi (.)\right|}{\left|\xi (.)\right|} \leq \epsilon $} 
\STATE Solve RMP to obtain final consumption schedule $\boldsymbol{p}_i^*,\ \forall i$.
\end{algorithmic}
\end{algorithm}

\section{Experiments}

In this study, we propose not only a DR strategy to control the load consumption of a residential community by formulating a MILP problem but also an efficient distributed optimization framework that can solve it even if the number of homes in the community is considerably large. Hence, we test those hypotheses separately by conducting
different set of experiments. Section \ref{exp_sim} provides the details of our
experimental setup; Section \ref{exp_load} shows that the proposed formulation is
useful for controlling the load consumption of the community; and Section
\ref{exp_comp} compares the optimization performance of the proposed algorithm with a
commercial solver.

\subsection {Simulation Setup}
\label{exp_sim}

We designed (in Python) a simulation environment that can create a community of
homes, each equipped with the appliances introduced in Appendix \ref{app_homeappliance}. We
utilize the Gurobi 9.5.1 package to solve the formulated MILP and the {\em Linear
  Programming (LP)} problems in Section \ref{section:opt}. Although the properties of
the appliances are identical across the community, the preferences of households for
each appliance can vary. Therefore, the preferences such as comfortable room
temperature and the desired water temperature are randomly initialized by using
different probability distributions. Moreover, the desirable load consumption
schedule, $\overline{\boldsymbol{p}_i}$, for each home is initialized in line with
the specified preferences and consumption habits of the corresponding home. Both the details of the appliance properties and the probability distributions that we use to sample the preferences of households are provided in Appendix \ref{app_homeappliance}. We randomly sample the rate of incentive payments, $c_{ij}$, for each home independently from a Normal Distribution with mean 0.01 and standard deviation 0.005. \textcolor{black}{We set the rate to 0 if the randomly sampled value is negative.} Similarly, we sample the aggregate uncontrollable demand of the community from a Normal Distribution with mean 300 and standard deviation 20 for each time point independently. We assume that PVs are renewable resources that provide energy to the community, so the power from renewable sources is available only from 6:00 to 18:00. %Accordingly, for this time period, we characterize the load generation of each PV at time $T$ by the polynomial $p^r(T) = -16T^2 + 384T - 1728$. 
Finally, the source code is publicly available\footnote{Code available at \url{https://github.com/erhancanozcan/energy-MILP}.}.%\footnote{URL will be included in the camera-ready version.}.

%For a given time index, we characterize the load generation of each PV at time t 
%and characterize the load generation of each PV at time index $t$ by the polynomial $p^r(t) = -t^2 + 96t - 1728$
%$p^r(t) = -16t^2 + 384t - 1728$

In our simulations, we assume that the length of the planning horizon is 24 hours,
and each time interval lasts 15 minutes $(K=96)$. Since the desirable load
consumption schedule of each home can change dynamically at each time interval, a practical
optimization algorithm should provide a reasonably good solution for each home in less
than 15 minutes, \textcolor{black}{which aligns with industry practices and guidelines \cite{qdr2006benefits}.}

%Lastly, the primary role of the proposed formulation is to maintain a balance in load demand and supply by purchasing additional power from external sources, and the performance of the algorithm can be subject  The amount of uncontrollable load and the generated power from PVs

% Finally, the role of the aggregator is to maintain a balance in load demand and supply by purchasing additional power from external sources.

% we consider a residential neighborhood whose
% power demand is met by both renewable resources and external
% traditional generators, and design a DR strategy to match sup-
% ply with demand while satisfying the user demand. 

% keep the total power consumption of the
% community close to $Q(t)$ at time interval $t$, and this value must be determined by
% the aggregator before optimizing the problem. In our simulations, we set
% $Q(t)$ as the time average of the total desirable power consumption to prevent
% instantaneous fluctuations, and $Q(t)$ takes the form
% \[
% Q(t)=\frac{\sum_{i=1}^N\sum_{j=1}^M\sum_{t=0}^{K-1}\overline{p_{ij}}(t)}{K},\ \forall
% t.
% \]
% However, in practice, the aggregator can determine the targeted aggregate
% power consumption level based on the amount of power supplied by the generators.

\subsection {The Load Shaping Capacity}
\label{exp_load}

In this section, we show the benefit of utilizing the proposed formulation to shape
the overall load consumption of the community, and analyze the performance of the algorithm under different scenarios by varying the amount of total power generated by renewables. When no optimization algorithm is
employed, each home simply follows its desirable load schedule. Therefore, the total load consumption of home $i$ at time $t$ is equal to
$\sum_{j=1}^M\overline{p_{ij}}(t) + p_{i,\text{U}}(t)$, and the desirable load consumption of the whole
community is equal to $\sum_{i=1}^N\sum_{j=1}^M\overline{p_{ij}}(t) + \sum_{i=1}^Np_{i,\text{U}}(t)$ when there is no optimization. On the other hand, $\sum_{i=1}^N\sum_{j=1}^M p_{ij}^*(t) + \sum_{i=1}^Np_{i,\text{U}}(t)$ denotes the optimal load consumption of
the community at time $t$ after optimizing the proposed formulation. Figure
\ref{fig:load_demand} shows the change in the load consumption of the community after solving the proposed optimization problem for two different PV generation scenarios. 
\begin{figure}[htb]
\centering
\includesvg[width=3.75in]{load_shaping_demand.svg}
\vspace{-1.5cm}
\caption{Change in total load consumption of the community with 1000 homes after optimization under low (top) and high (bottom) PV generation scenarios.} 
\label{fig:load_demand}
\end{figure}

According to Figure \ref{fig:load_demand}, we can observe that the overall load demand of the community changes significantly after optimization for both scenarios. We observe an increase in the load consumption of the community around noon under the high renewable generation scenario. However, this is not an indication of peak problem, since the increase in load demand can be satisfied by the renewable resources. In our formulation, the goal of the aggregator is to keep the amount of purchased power from external sources constant during the planning horizon. Therefore, in order to understand the functionality of the proposed formulation, we calculate the difference between total demand and renewable generation before and after optimization, and compare these in Figure \ref{fig:load_demand_optimalQ}.
\begin{figure}[!h]
\centering
\includesvg[width=3.75in]{load_shaping_optimalQ.svg}
\vspace{-1.5cm}
\caption{Change in the difference between total demand and renewables for a community with 1000 homes after optimization.} 
\label{fig:load_demand_optimalQ}
\end{figure}

According to Figure \ref{fig:load_demand_optimalQ}, for both low and high renewable generation scenarios, the difference between total demand and renewables remains close to the calculated $q^*$, e.g., $\sum_{i=1}^N\sum_{j=1}^M p_{ij}^*(t) + \sum_{i=1}^Np_{i,\text{U}}(t) - p^{r}(t) \sim q^*$, after optimization. On the other hand, when the participating homes determine their own power consumption schedules based on
their preferences without receiving any signal, the amount of the power aggregator purchases from external generators varies significantly, and the peak load problem remains untreated. 

Finally, since the load shaping capacity of the proposed formulation depends on the rate of incentive payment coefficients, $c_{ij}$, we explore the sensitivity of the proposed formulation to the rate of incentive payment coefficients. To measure the load shaping capacity of the proposed formulation, we calculate the mean absolute deviation between total demand and total supply across time points, e.g., {$\text{MAD}=\nobreak \frac{\sum_{t=0}^{K-1}\left|\left(\sum_{i=1}^N\sum_{j=1}^M p_{ij}^*(t) + \sum_{i=1}^Np_{i,\text{U}}(t)\right)- \left(p^{r}(t) + Q^*(t)\right)\right| }{K}$}, and observe the change of MAD by varying incentive payment rates. In our experiments, we sample these rates from a Normal Distribution for each home independently. Thus, we consider different mean values for incentive payment rates. Finally, to decrease the effect of random initializations, we repeat each experiment for five different seeds, and summarize our findings via box plots in Figure \ref{fig:dev_cost_sensitivity}. According to Figure \ref{fig:dev_cost_sensitivity}, as the incentive payment rates increase, minimizing the distributed incentive payments becomes more important for the optimization. As a result, we start observing an increase in the MAD values. On the other hand, when the incentive payment rates are small, maintaining the balance between supply and demand becomes easier for the aggregator.
% \begin{equation}
% \text{MAD}= \frac{\sum_{t=0}^{K-1}\left|\left(\sum_{i=1}^N\sum_{j=1}^M p_{ij}^*(t) + \sum_{i=1}^Np_{i,\text{U}}(t)\right)- \left(p^{r}(t) + Q^*(t)\right)\right| }{K} 
% \end{equation}
\begin{figure}[!h]
\centering
\includesvg[width=3.25in]{dev_cost_sensitivity.svg}
\vspace{-0.5cm}
\caption{Load shaping capacity of the formulation under different rates of incentive payments.} 
\label{fig:dev_cost_sensitivity}
\end{figure}

\subsection {Computational Performance of Distributed Optimization}
\label{exp_comp}

Despite the effectiveness of the proposed formulation to control the load consumption
of the community, the proposed formulation is a MILP problem, and it is hard to solve
this problem when the number of homes in the community is large. Since the main goal
of this study is to maintain a balance between load supply and demand in residential neighborhoods, the
number of participating homes in large neighborhoods can be on the scale of
thousands. Hence, it is essential to design an efficient algorithm to solve the proposed optimization problem.
% Furthermore, the load consumption plans of the households change
% dynamically. In order to accommodate those changes, we propose re-optimizing the
% problem before the next time interval starts. Hence, an efficient algorithm is
% required to solve the proposed optimization problem.

Gurobi can optimize the formulated Centralized-IP \textcolor{black}{in \eqref{central_ip}} when the community size is
small. However, the optimization time significantly increases for large
communities. The optimization duration can be shortened by selecting larger
optimality gaps for Gurobi. The default value of this parameter in Gurobi is
$1e^{-4}$, and we set it to $1e^{-2}$ in our experiments. On the other hand, the parameter $\epsilon$ playing a similar role
to the optimality gap of Gurobi is set to $1e^{-3}$ in our algorithm. In addition to
this, our algorithm has a hyperparameter $\kappa$ controlling the column removal
process, and the computational performance of our algorithm can be improved
substantially by selecting this parameter carefully. By setting $\kappa$ to a large
number, we obtain a naive version of our algorithm, which does not remove any column
from the problem. In addition to this, we obtain two more versions of our algorithm
by setting $\kappa$ equal to $5$, and $10$, respectively. In order to compare the
optimization time of the alternative approaches, we randomly initialize communities
with varying number of homes (e.g., $1000$, $3000$, $5000$, and $7000$). To decrease
the effect of random initialization, we use $5$ different seeds. Figure
\ref{fig_opt_time} shows the mean optimization time (along with a $95$\% confidence
interval) for each algorithm.

\begin{figure}[htb]
\centering
\includesvg[width=3.5in]{time_complexity.svg}
\caption{Change in optimization time with respect to community size.}
\label{fig_opt_time}
\end{figure}

The load consumption plans of the households change dynamically. In order to accommodate those changes, it might be necessary to re-optimize the
problem before the next time interval starts for a practical optimization algorithm while keeping $Q^*$ fixed. Therefore, in our case, the run time of the algorithm has to be less than $15$ minutes ($900$ seconds). According to Figure \ref{fig_opt_time}, Gurobi fails to satisfy this requirement. On the other hand, the computational complexity of our
distributed optimization algorithm increases linearly, and it can provide a solution
within the time limit when unused columns are removed from the relaxed-RMP. At this
point, we want to make a comment on the optimization time of the distributed
algorithm, which does not remove any of the columns. In Section \ref{section:opt}, it
is argued that removing unused columns from the relaxed-RMP helps our algorithm to
solve the relaxed-RMP within a reasonably short time. However, the last step of our
algorithm is to solve the RMP with the available columns, and this operation can be
costly when the number of columns is excessive. Removing unused columns during the
optimization of the relaxed-RMP enables us to formulate a simpler RMP in the last
step of our algorithm. Thus, this step also takes less time, and the overall
optimization time of our distributed algorithm decreases significantly.

%Hence, an efficient algorithm is required to solve the proposed optimization problem.

The quality of the solutions provided by each algorithm is important as well as the
optimization time of algorithms. Therefore, we have to evaluate the quality of the
solutions provided by algorithms. Suppose that we have an oracle that can optimally
solve the given problem, and the optimal objective value for this problem is denoted
by $z^*$. In practice, we solve the problem by using Gurobi to obtain this value. On the other hand, we
have mentioned 3 different versions of our algorithm to optimize the problem in this section, and
$z^{\text{method}_k}$ represents the best objective value that the k-th method can
find. Then, we define the pseudogap of the k-th method as follows:
\[
\text{pseudogap (}\text{method}_k) =
\frac{z^{\text{method}_k}-z^*}{z^*},\ \forall\ k.
\]
Figure \ref{fig_opt_gap} shows the mean pseudogap (along with a 95\% confidence
interval) for each algorithm.

\begin{figure}[htb]
\centering
\includesvg[width=3.5in]{optimality_gap_nogurobi.svg}
\caption{Change in the pseudogap with respect to community size.}
\label{fig_opt_gap}
\end{figure}

Based on Figure \ref{fig_opt_gap}, our distributed heuristic optimization framework provides the smallest optimality gap when all of the columns are maintained, which is expected as the algorithm has more flexibility to solve the RMP in the last step. However, similar to Gurobi, this version of our algorithm is not practical as it cannot provide the solution within the desired time limit. On the other hand,  our distributed heuristic optimization framework can efficiently solve the
problem even if the community size is huge by utilizing the column removal idea, and the pseudogaps of these algorithms are always less than 0.01 according to Figure \ref{fig_opt_gap}. Therefore, we can conclude that the proposed distributed
optimization framework can provide reasonably good solutions in significantly less
amount of time.

\section{Conclusions}

We presented a new MILP formulation that can control the load consumption schedules of
households in a neighborhood while integrating renewable resources. Since the load consumption of each appliance is
determined according to the preferences of individuals, we guarantee that the user
comfort is always maintained. We developed a distributed approach based on a
Dantzig-Wolfe decomposition approach to solve the formulated problem efficiently in
large communities. In addition to improving the optimization time, our distributed
optimization framework helps us to address the possible concerns related to data
privacy. In the future, our goal is to investigate the effect of using a quadratic
penalty function to minimize the deviations from the targeted aggregate power
consumption level. Then, we will test the performance of our distributed optimization
framework on this new version. Finally, we also want to formulate a robust problem considering uncertainties in renewable resources to make our approach more realistic.

\appendices
\section{Home Appliances}
\label{app_homeappliance}
The set $\H_i$ consists of the comfort related constraints of home $i$, and we
introduce the necessary constraints to maintain the user comfort for each appliance
in this section. The planning horizon is divided into $K$ equal time intervals, and
the whole planning horizon is spanned by the set $\T:=\{0,1,\dots,K-1\}$.
% where
% \[
% \T:=\{0,1,\dots,K-1\}.
% \]

\subsection{HVAC}
We assume that each home has a built-in HVAC system,which operates according to simple logical rules. Suppose that consumer $i$ feels comfortable as long as the room temperature is
between $T_i^{low}$ and $T_i^{upper}$. Hence, whenever the HVAC is in its heating
mode and the room temperature is less than the lower temperature limit, the HVAC gets
activated. On the other hand, the HVAC starts working if the room temperature is
greater than the upper temperature limit and it is in its cooling mode. Therefore,
designing a mechanism, which keeps the room temperature in the comfortable thermal
range is our main goal. 

While designing the HVAC model, we simplify the temperature update rule given in \cite{foresee2017} by
omitting the solar irradiance effect, and use the following
equation to update the indoor room temperature:

\begin{equation*}
    T_{in}(t+1)=T_{in}(t) + \gamma_1\left(T_{out}(t)-T_{in}(t)\right) + \gamma_2 p_{\text{HVAC}}(t) m,
\end{equation*}
where $T_{in}(t)$ is the room temperature at time $t$, $T_{out}(t)$ is the outside temperature at time $t$, $p_{\text{HVAC}}(t)$ is the consumed power at time t, $m$ represents either the heating (1) or the cooling (-1) mode of HVAC, and $\gamma_1, \gamma_2$ are the home specific coefficients
related to building isolation and heating or cooling gain, respectively.

However, it is hard to ensure that the installed HVAC will
always satisfy the temperature requirements of the home since there are various
factors affecting the thermal dynamics of a room. Also, in some extreme weather
conditions, the HVAC might be inadequate either to warm-up or to cool down the room
immediately since its capacity is limited. Due to aforementioned reasons, we focus on an approach limiting the
deviations from the comfortable thermal range.

Let us assume that consumer $i$ controls the HVAC by only defining her comfortable
thermal range without using any optimization algorithm. Then, the expected room
temperature, $T_{i,in}^e(t)$, can be calculated recursively by using the following
equations:
\begin{align*}
 T_{i,in}^e(t+1) = & T_{i,in}^e(t) +\gamma_1\left(T_{out}(t)-T_{i,in}^e(t)\right) \\
 &\ + \gamma_2 p_{i,\text{HVAC}}(t),\qquad t \in \T,\\ 
 p_{i,\text{HVAC}}(t) = &a_{i,\text{HVAC}}(t) \alpha_{\text{HVAC}}  m
  N^{m}_{\text{HVAC}}, \qquad t \in \T,
\end{align*}
where
\begin{equation*}
a_{i,\text{HVAC}}(t)=\begin{cases}
                        1,&{\text{if}}\ T_{i,in}^e(t)\leq T_i^{low}\ \&\ m=1 \\ 
                        1,&{\text{if}}\ T_{i,in}^e(t)\geq T_i^{upper}\ \&\ m=-1 \\ 
                        {0,}&{\text{otherwise}}
\end{cases},
\end{equation*}
$\alpha_{\text{HVAC}}$ is thermal conversion efficiency of HVAC, $m$ is the heating (1) or cooling (-1) mode of HVAC, and $N_{\text{HVAC}}^{m}$ is the nominal power consumption of HVAC when the mode $m$ is active.
%\item[$\alpha_{\text{HVAC}}$] Thermal conversion efficiency of HVAC
%\item[$m$] Heating (1) / Cooling (-1) mode of HVAC
%\item[$N_{\text{HVAC}}^{m}$] Nominal power consumption of HVAC when the mode $m$ is active
By using the expected room temperature values, deviations from the comfortable range
are calculated as follows:
%\label{temp_dev}
\begin{align*}
&s_{i,+}^*(t) = max\left(T_{i,in}^e(t) - T_i^{upper}),0\right),\ t \in \T,\\
&s_{i,-}^*(t) = max\left(T_i^{low} -T_{i,in}^e(t) ,0\right),\ t \in \T.
\end{align*}
Each home $i$ has parameter, $\epsilon_{i,HVAC}$, which lets deviations from the
comfortable thermal range to be slightly larger. This parameter gives some
flexibility to the HEMS while setting the working schedule of HVAC. We define the
HVAC related decision variables as follows:
\begin{itemize}
    \item $a_{i,\text{HVAC}}(t)$ denotes the ON/OFF status of HVAC at time $t$.
    \item $s_{i,-}(t))$ denotes the amount of deviation from the lower temperature limit at time $t$.
    \item $s_{i,+}(t))$ denotes the amount of deviation from the upper temperature limit at time $t$.
    \item $T_i^{set}(t)$ denotes the set temperature HVAC needs to operate.
\end{itemize}
Then, HEMS can determine the working schedule of HVAC in the heating mode by
considering the following constraints:
\begin{subequations}
\begin{align}
  \notag \quad& T_{i,in}(t+1)=T_{i,in}(t) +\gamma_1(T_{out}(t)-T_{i,in}(t)) \\
 &\ \quad\quad\quad\quad\quad
  +\alpha_{\text{HVAC}}m\gamma_2p_{i,\text{HVAC}}(t),\quad t\in \T,\label{hvac_b}\\ 
\quad&p_{i,\text{HVAC}}(t)= a_{i,\text{HVAC}}(t) N^{m}_\text{HVAC},\quad t \in
\T,\label{hvac_c} \\ 
\quad& T_i^{low} - s_{i,-}(t) \leq T_{i,in}(t),\quad t\in \T,\label{hvac_d} \\
\quad& T_{i,in}(t)\leq T_i^{upper} + s_{i,+}(t),\quad t\in \T,\label{hvac_e}\\
\quad& s_{i,-}(t) \leq s_{i,-}^*(t) + \epsilon_{i,\text{HVAC}},\quad t\in \T,\label{hvac_f}\\
\quad& s_{i,+}(t) \leq s_{i,+}^*(t) + \epsilon_{i,\text{HVAC}},\quad t\in \T,\label{hvac_g}\\
\quad& T_i^{set}(t) \geq T_{i,in}(t) - M^{H}(1-a_{i,\text{HVAC}}(t)),\quad t\in
\T,\label{hvac_h}\\ 
\quad& T_{i,in}(t) \geq T_i^{set}(t)  - M^{H} a_{i,\text{HVAC}}(t),\quad t\in
\T,\label{hvac_i}\\ 
\quad& s_{i,-}(t),s_{i,+}(t)\geq0,\quad t\in \T,\label{hvac_j}\\
\quad& a_{i,\text{HVAC}}(t) \in \{0,1\},\quad t\in \T\label{hvac_k},
\end{align}
\end{subequations}
where Equation \eqref{hvac_b} estimates the room temperature recursively based on
other parameters, Equation \eqref{hvac_c} calculates the power consumed by HVAC at
time $t$, and Equations \eqref{hvac_d}, \eqref{hvac_e}, \eqref{hvac_f}, and
\eqref{hvac_g} ensure that the room temperature does not deviate too much from the
comfortable thermal range at time $t$. Finally, Equations \eqref{hvac_h},
\eqref{hvac_i}, and \eqref{hvac_k} relate the set temperature and ON/OFF status of
HVAC with the help of some large number $M^H$. Although Equation \eqref{hvac_k}
states that $a_{i,\text{HVAC}}(t)$ is a binary variable, we relaxed this constraint
in our implementation by letting $a_{i,\text{HVAC}}(t)$ to be a non-negative real
number. This is a reasonable relaxation because we can also interpret this variable
as the percentage of time on which the HVAC remains active at time interval $t$. Even
without this relaxation, our optimization framework would work with the
constraint set given in Equations \eqref{hvac_b}-\eqref{hvac_k} because the set
$\X_i$ we consider is already a mixed integer polyhedron.

While using the HVAC in its cooling mode, the constraints in Equations \eqref{hvac_h}
and \eqref{hvac_i} must be replaced by the following two equations:
\begin{align*}
&T_{i,in}(t) \geq T_i^{set}(t) -M^{H}(1-a_{i,\text{HVAC}}(t)),\\
&T_{i,in}(t) \leq T_i^{set}(t)  +M^{H}a_{i,\text{HVAC}}(t).
\end{align*}

Finally, Table \ref{tab:HVAC} summarizes the properties of the HVAC and random
distributions that we utilized to have a diversified neighborhood.
\begin{table}[h]
\caption{The Properties of HVAC\label{tab:HVAC}}
\centering
\renewcommand{\arraystretch}{1.5}
\begin{tabular}{|c|c|}
\hline
$\gamma_1$ & $N\left(0.10,0.001\right)$ \\
\hline
$\gamma_2$ & $N\left(3e^{-6},1e^{-7}\right)$\\
\hline
$T_{out}$ & Weather Data from 2018 \cite{my_ref18} \\
\hline
$\alpha_{i,\text{HVAC}}$ & $0.9$ \\
\hline
$N^m_{\text{HVAC}}$ & $m=1$ $\Rightarrow$ $3$ (kW), $m=-1$ $\Rightarrow$ $2$ (kW) \\
\hline
$T_i^{low}$ & Discrete U[19,24] $(^\circ C)$\\
\hline
$T_i^{upper}$ & $T_i^{low}+2$ $(^\circ C)$\\
\hline
$\epsilon_{i,\text{HVAC}}$ & $0.5$ $(^\circ C)$\\
\hline
\end{tabular}
\end{table}

\subsection{EWH}

Similar to \cite{my_ref13}, we assume that the power consumption of the EWH on
standby is negligible because the water tanks are well insulated. On the other hand,
the EWH consumes significant amount of power while heating the tap water until the
desired temperature specified by the user. Unlike the EWH model introduced in
\cite{my_ref13}, we model the amount of available hot water in the tank instead of
the average water temperature in the tank. To that end, we assume that
there are two separate tanks, and the water in those tanks does not mix unless the
water temperatures in both tanks is the same. In other words, while one
of the tanks serves as the main storage tank, the heating occurs in the other tank,
and the tap water is transferred to the main storage tank without any heat loss when
it reaches the desired temperature. This additional assumption allows us to update
the amount of available hot water in the main tank at each time interval by using the
relation between heat transfer and temperature change. Suppose that
$z_{\text{EWH}}(t)$ denotes the amount of water that the EWH system can heat from
time $t$ to time $(t+1)$. Then, the relation between $p_{\text{EWH}}(t)$ and
$z_{\text{EWH}}(t)$ is as follows:
\begin{equation}
z_{\text{EWH}}(t) = \frac{p_{\text{EWH}}(t)\eta_{\text{EWH}}}{\rho(T_{d}-T_{t})},
\end{equation}
where $\eta_{\text{EWH}}$ is water heating efficiency of EWH, $\rho$ is specific heat of water $\left(J / kg^\circ C\right)$, $T_{d}$ is desired water temperature $(^\circ C)$, and $T_{t}$ is tap water temperature $(^\circ C)$.

Suppose that $C_{i,\text{EWH}}$ denotes the capacity of the tank,
$p_{\text{EHW}}^{max_i}$ denotes the maximum power the EWH can consume, and
$y_{i,EWH}(t)$ represents the hot water demand of home $i$ at time point $t$. We
define the following decision variables:
\begin{itemize}
    \item $x_{i,\text{EWH}}(t)$ denotes the amount of hot water available in tank at time $t$.
    \item $z_{i,\text{EWH}}(t)$ denotes the amount of water that the EWH will heat from time $t$ to $(t+1)$.
    \item $p_{i,\text{EWH}}(t)$ denotes the amount of power provided by the EWH at time $t$.
\end{itemize}

Then, the hot water consumption preferences of the user $i$ can be satisfied by adding the following set of constraints to the HEMS:
\begin{subequations}
\begin{align}
&x_{i,\text{EWH}}(t) \leq C_{i,\text{EWH}}, \quad t \in \T, \label{ewh_a}\\
&x_{i,\text{EWH}}(t) \geq y_{i,\text{EWH}}(t), \quad t \in \T,\label{ewh_b}\\
&p_{i,\text{EWH}}(t) \leq p_{\text{EHW}}^{max_i},\quad t \in \T,\label{ewh_c}\\
&p_{i,\text{EWH}}(t) \geq 0,\quad t \in \T,\label{ewh_d}\\
&z_{i,\text{EWH}}(t) = \frac{p_{i,\text{EWH}}(t)\; \eta_{\text{EWH}}}{\rho(T_{d} - T_{t})}, \quad t \in \T,\label{ewh_e}\\
&x_{i,\text{EWH}}(t) \geq 0,\quad t \in \T,\label{ewh_f}\\
\notag &x_{i,\text{EWH}}(t+1) = x_{i,\text{EWH}}(t) + z_{i,\text{EWH}}(t)\\ &\ \quad\qquad\qquad\qquad\qquad - y_{i,EWH}(t),  t \in \T,\label{ewh_g}
\end{align}
\end{subequations}
where Equation \eqref{ewh_a} ensures that the available hot water does not exceed the
capacity of the tank, Equation \eqref{ewh_b} ensures that the water tank has enough
hot water to satisfy the water demand, Equations \eqref{ewh_c} and \eqref{ewh_d} set
the lower and upper bounds for the power consumption of the EWH, Equation
\eqref{ewh_e} calculates the amount of tap water that the EWH can heat, and Equation
\eqref{ewh_f} ensures that the amount of water in the tank is always
non-negative. Finally, Equation \eqref{ewh_g} defines the amount of available hot
water in the tank recursively by considering the inflow and outflow of the hot
water. The parameters of the EWH are summarized in Table \ref{tab:EWH}.

\begin{table}[h]
\caption{The Properties of EWH\label{tab:EWH}}
\centering
\renewcommand{\arraystretch}{1.5}
\begin{tabular}{|c|c|}
\hline
$C_{i,\text{EWH}}$ & $270$ (kg) \\
\hline
$y_{i,\text{EWH}}(t)$ & $| N\left(30,10\right) |$ (kg) $, \forall\ t \in \T^{\text{EWH}_i}$\\
\hline
$\T^{\text{EWH}_i}$ & Consists of $n_{i,EWH}$ random points selected from $\T$\\
\hline
$n_{i,\text{EWH}}$ & Discrete U[2,5]\\
\hline
$T_{d}$ & Discrete U[40,42] $(^\circ C)$ \\
\hline
$T_{t}$ & $4$ $(^\circ C)$ \\
\hline
$\eta_{\text{EWH}}$ & $0.95$ \\
\hline
$\rho$ & $4186 \left(J / kg^\circ C\right)$ \\
\hline
$p_{\text{EHW}}^{max_i}$ & $4$ (kW) \\
\hline
\end{tabular}
\end{table}

\subsection{EV}

We assume that the consumer is flexible enough to charge the battery of their vehicle
whenever the vehicle is not in use, and has an access to a 240-volt outlet to charge
the battery. To prevent any possible equipment degradation issue in the long run, the
consumer can limit the amount of instantaneous maximum amperage delivered to the
battery. HEMS determines the amount of the current that will be delivered to the
battery by considering the battery {\em State Of Charge (SOC)}, required power to
complete the scheduled next trips, and the maximum charging rate. Finally, due to
Ohm's Law, the relation between the current and the power delivered to the battery is
as follows:
\begin{equation}
p_\text{EV}(t)=\frac{240V\times I(t)}{1000},
\end{equation}
where $I(t)$ is the current delivered to the EV battery at time $t$ $(A)$.

Suppose that $C_{i,\text{EV}}$ denotes the power capacity of the EV battery,
$I^{max_i}$ denotes the maximum amperage that can be delivered to the battery,
$y_{i,\text{EV}}(t)$ denotes the amount of power EV needs from time point $t$ to
$(t+1)$. Finally, since the battery of EV cannot be charged when it is in use, we
have to define the set $\T^{\text{EV}_i}:=\{t\ |\ y_{i,\text{EV}}(t)>0,\ t\in \T\}$.
We define the decision variables for the EV as follows:
\begin{itemize}
    \item $x_{i,\text{EV}}(t)$ denotes the amount of available power in the battery at time $t$.
    \item $I_i(t)$ denotes the amount of current delivered to the battery of the EV from time $t$ to $(t+1)$.
\end{itemize}
Then, HEMS can satisfy the EV related requirements of home $i$ when the followings set of constraints are added to the HEMS:
\begin{subequations}
\begin{align}
& x_{i,\text{EV}}(t) \leq C_{i,\text{EV}}, \quad t \in \T, \label{ev_a}\\
& x_{i,\text{EV}}(t) \geq y_{i,\text{EV}}(t), \quad t \in \T, \label{ev_b}\\
& p_{i,\text{EV}}(t) = \frac{240\text{Volt}*I_i(t)}{1000},\quad t \in \T, \label{ev_c}\\
& p_{i,\text{EV}}(t) \geq 0, \quad t \in \T, \label{ev_d}\\
& I_i(t) \leq I^{max_i}, \quad t \in \T, \label{ev_e}\\
& I_i(t) = 0,\quad t \in T^{\text{EV}_i}, \label{ev_f}\\
& x_{i,\text{EV}}(t+1) = x_{i,\text{EV}}(t) + p_{i,\text{EV}}(t) - y_{i,\text{EV}}(t) \quad t \in \T, \label{ev_g}
\end{align}
\end{subequations}
where Equation \eqref{ev_a} ensures that the amount of power stored in the battery
does not exceed its capacity, Equation \eqref{ev_b} forces that enough power is
stored in the battery to satisfy the customer needs, Equation \eqref{ev_c} calculates
the amount of power delivered to the battery at time $t$, Equation \eqref{ev_e}
limits the amount of maximum amperage delivered to the battery, Equation \eqref{ev_f}
guarantees that the battery cannot be charged when the EV is in use. Finally,
Equation \eqref{ev_g} defines the amount of stored power in the battery recursively
by considering charging and discharging events. The parameters of the EV are
summarized in Table \ref{tab:EV}.

\begin{table}[!t]
\caption{The Properties of EV\label{tab:EV}}
\centering
\renewcommand{\arraystretch}{1.5}
\begin{tabular}{|c|c|}
\hline
$C_{i,\text{EV}}$ & $60$ (kWh) \\
\hline
$y_{i,\text{EV}}(t)$ & $0.346\times x_{i,\text{EV}}(t)$ (kWh) $, \forall\ t \in \T^{\text{EV}_i}$\\
\hline
$x_{i,\text{EV}}(t)$ & Discrete U[5,9] (miles) $, \forall\ t \in \T^{\text{EV}_i}$\\
\hline
$\T^{\text{EV}_i}$ & Consists of $n_{i,\text{EV}}$ random points selected from $\T$\\
\hline
$n_{i,\text{EV}}$ & Discrete U[4,12]\\
\hline
$I^{max_i}$ & $24$ $(A)$ \\
\hline
\end{tabular}
\end{table}

\subsection{Basic Appliances}

Basic appliances such as oven, washing machine, and dryer offer flexible
operation. In other words, the operation of these appliances can be either scheduled
to start early or delayed as long as it is completed in the user specified time
interval. We model these appliances similarly to the delayed flexible appliances in
\cite{chen2013mpc}, so that the appliance runs without any interruption once it
starts. However, since the continuity of the process can only be ensured by using
binary variables in the constraints, the feasible region defining the customer
comfort becomes non-convex.

Basic appliances that we consider in this study are: Washing Machine (WM), oven, and
dryer, and they are all similarly modeled. In this section, we define the constraints
of the WM, but these constraints can be easily generalized for other
appliances. Basic appliances have ON / OFF status, and we assume that the power
consumption is constant during the operation.

Suppose that home $i$ wants WM operation to be completed between
$\alpha^{start}_{i,\text{WM}}$ and $\alpha^{end}_{i,\text{WM}}$. Accordingly, we
define the set $\T^{\text{WM}_i}:=\{t\ |\ \alpha^{start}_{i,\text{WM}}\leq t \leq
\alpha^{end}_{i,\text{WM}},\ t\in \T\}$. Moreover, let $N_{i,\text{WM}}$ represents
the nominal power consumed by the WM when it is operating, and the operation lasts
$k_{WM}$ periods without any interruption once it starts. In order to enforce the
continuity of the process, we need to define the following binary variables:
\begin{itemize}
    \item $x_{i,\text{WM}}(t)$ denotes ON $(1)$ / OFF $(0)$ status at time $t$.
    \item $z_{i,\text{WM}}(t)$ denotes the start-up action at the beginning of time period $t$ when it is equal to $1$.
    \item $y_{i,\text{WM}}(t)$ denotes the shut-down action at the beginning of time period $t$ when it is equal to $1$.
\end{itemize}
We add the following set of constraints to the HEMS to ensure that WM operates according to the preferences of home $i$.
\begin{subequations}
\begin{align}
& \sum_{t \in \T^{\text{WM}_i}} p_{i,\text{WM}}(t)= N_{i,\text{WM}}\ .\ k_{\text{WM}}, \label{wm_a}\\
& p_{i,\text{WM}}(t) = N_{i,\text{WM}}\ .\ x_{i,\text{WM}}(t),\quad t \in \T^{\text{WM}_i},\label{wm_b}\\
& x_{i,\text{WM}}(t) = 0,\quad t \in \T\setminus\T^{\text{WM}_i},\label{wm_c}\\
\notag & x_{i,\text{WM}}(t) - x_{i,\text{WM}}(t-1) = z_{i,\text{WM}}(t) \\ 
\notag                                    &\ \qquad\qquad\qquad\qquad\qquad\quad- y_{i,\text{WM}}(t),\\ 
                                        &\ \qquad\qquad\qquad\qquad\qquad\quad t \in \T^{\text{WM}_i}\cup\{\alpha^{end}_{i,\text{WM}}+1\},\label{wm_d}\\
& z_{i,\text{WM}}(t) \leq x_{i,\text{WM}}(t), \quad t \in \T^{\text{WM}_i},\label{wm_e}\\
& y_{i,\text{WM}}(t)+x_{i,\text{WM}}(t) \leq 1, \quad t \in \T^{\text{WM}_i},\label{wm_f}\\
& \sum_{t \in \T^{\text{WM}_i}} z_{i,\text{WM}}(t) \leq 1,\quad\label{wm_g}\\
& \sum_{t \in \T^{\text{WM}_i}} y_{i,\text{WM}}(t) \leq 1,\quad\label{wm_h}\\
 &x_{i,\text{WM}}(t),y_{i,\text{WM}}(t),z_{i,\text{WM}}(t) \ \  \in \{0,1\},\ \forall \ t \in \T^{\text{WM}_i},\label{wm_i}
\end{align}
\end{subequations}
where Equation \eqref{wm_a} ensures that the WM operates for $k_{\text{WM}}$ periods,
Equation \eqref{wm_b} calculates the power consumed by the WM at time period $t$,
Equation \eqref{wm_c} prevents the WM from running when $t \notin \T^{\text{WM}_i}$,
Equation \eqref{wm_d} sets the logical relations between start-up/shut-down and
ON/OFF status, Equation \eqref{wm_e} ensures that the WM is in the ON status whenever
a start-up action takes place, and Equation \eqref{wm_f} ensures that the WM is in
the OFF status if a shut-down action takes place in the beginning of the
period. Equations \eqref{wm_g} and \eqref{wm_h} guarantee that start-up and shut-down
can occur at most once during the planning horizon. Finally, we define the required
binary variables in Equation \eqref{wm_i}. Table \ref{tab:WM} shows the parameters of
the WM.

\begin{table}[h]
\caption{The Properties of WM\label{tab:WM}}
\centering
\renewcommand{\arraystretch}{1.5}
\begin{tabular}{|c|c|}
\hline
$N_{i,\text{WM}}$ & $0.5$ (kWh) \\
\hline
$k_{\text{WM}}$ & 4 \\
\hline
$\alpha_{i,\text{WM}}^{start}$ & max$(s_{i,\text{WM}}-4,0)$ \\
\hline
$\alpha_{i,\text{WM}}^{end}$ & min$(s_{i,\text{WM}}+8, K-1)$ \\
\hline
$s_{i,\text{WM}}$ & Discrete U[0,K-1] \\
\hline
\end{tabular}
\end{table}

\section{Distributed Optimization}
\label{app_distopt}

The RMP defined in Equations \eqref{eq_13}-\eqref{rmp_last} is a binary integer programming
problem. We relax this problem by removing the integrality constraints on the
$\lambda_i^g$ variables, and the resulting relaxed-RMP is as follows:

\noindent\textbf{relaxed-RMP:}
\begin{subequations}
\begin{align}
\min_{\lambda_i^g, q, \boldsymbol{s}, \boldsymbol{a}} &\sum_{t=0}^{K-1} s(t) + \sum_{i=1}^N\sum_{g \in \widehat{G_i}}\left(\sum_{j=1}^M\sum_{t=0}^{K-1} c_{ij} u_{ij}^{+,g}(t)\right)\lambda_i^g\label{app_eq_13}  \\
\text{s.t.} \quad &s(t) - a(t) \geq \label{app_eq_rmp_first}
0,\quad\quad\quad\quad\quad\quad\quad\quad\quad\quad\quad \forall\ t, \\
&s(t) + a(t) \geq 0,\quad\quad\quad\quad\quad\quad\quad\quad\quad\quad\quad
\forall\ t, \\
&a(t)-q \notag\\
&\qquad+\sum_{i=1}^N\sum_{g \in \widehat{G_i}}\left(\sum_{j=1}^M
p_{ij}^g(t)\right)\lambda_i^g\notag\\
&\qquad=p^{r}(t) - \sum_{i=1}^Np_{i,\text{U}}(t),\ \qquad\qquad \forall\ t, \\
&\sum_{g \in
  \widehat{G_i}}\lambda_i^g=1,\quad\quad\quad\quad\quad\quad\qquad\quad\quad \forall\ i,\\
%&q(0) = q(1) = \cdots = q(K-1), \label{app_eq_rmp_last} \\
&q \geq 0, \qquad\qquad\qquad\qquad\quad\qquad\quad \\
&\lambda_i^g \geq 0,\qquad\qquad\qquad\qquad\qquad\quad \forall i, \forall g\in \widehat{G_i}.
\end{align}
\end{subequations}

Also, optimizing the relaxed-RMP by using the column generation mechanism described
in Section \ref{section:opt} is equivalent to optimizing the following convexified
version of Centralized-IP:

\noindent\textbf{Conv-IP:}
\begin{subequations}
\begin{flalign}
\min_{\boldsymbol{p}_1,\cdots,\boldsymbol{p}_N,\atop \boldsymbol{u}^+_1,\cdots,\boldsymbol{u}^+_N, q, \boldsymbol{s}, \boldsymbol{a}} &\sum_{t=0}^{K-1} s(t) + \sum_{i=1}^N \sum_{j=1}^M\sum_{t=0}^{K-1}c_{ij} u_{ij}^+(t)  \\[1em]
\text{s.t.} \quad &s(t) - a(t) \geq 0,\quad\quad\quad\quad\quad\quad\ \forall\ t\\
&s(t) + a(t) \geq 0,\quad\quad\quad\quad\quad\quad\ \forall\ t\\
&a(t)+\sum_{i=1}^N\sum_{j=1}^M p_{ij}(t)=q\ + \notag\\
&\quad\quad\quad\quad p^{r}(t) - \sum_{i=1}^Np_{i,\text{U}}(t),\quad \forall\ t,\\
%&q(0) = q(1) = \cdots = q(K-1),\\
&q \geq 0, \qquad\qquad\qquad\qquad\quad \\
&\left(\boldsymbol{p}_i,\boldsymbol{u}^+_i\right) \in \text{conv}(\X_i),\quad\quad\quad\ \ \ \forall\ i, \label{eq_39}
\end{flalign}
\end{subequations}
where conv$(\X_i)$ in Equation \eqref{eq_39} represents the convex hull of $\X_i$. 

In our case, since conv$(\X_i)$ is tighter than the feasible region expressed by the
LP relaxation of $\X_i$, the dual bounds that we calculate at each column generation
iteration are stronger. The following example shows why the region expressed by
conv$(\X_i)$ is tighter.

The comfort related constraints of all the appliances except basic appliances define
a convex region. Therefore, the convex hull and the LP relaxation can only be
different for basic appliances, and we design our example by focusing on the WM
related preferences of home $i$.\cmmnt{The constraints related to basic appliances
  make the proposed optimization problem harder since they involve binary variables,
  and we design our example by focusing on the WM related preferences of home $i$.}
The constraints in Equations \eqref{wm_a}-\eqref{wm_i} form $\H_i^{\text{WM}}$, and
home $i$ shares the feasible points in $\H_i^{\text{WM}}$ based on the signals sent
by the aggregator. Suppose that the length of the planning horizon is 5,
e.g., $\T:=\{0,1,\dots,4\}$, and the set $\T^{\text{WM}_i}$ includes the time points
when user $i$ can operate the WM, e.g.,$\T^{\text{WM}_i}:=\{1,2,3\}$. Also, let us
assume that the WM needs to work in two time periods to complete its cycle, and the
nominal power consumption of the WM, $N_{\text{WM}}$, is equal to 1.5. Based on these
assumptions, there are only two feasible points in $\H_i^{\text{WM}}$, and we provide
these points below: 
\begin{equation*}
A=\left\{
\begin{array}{lccc}
p(0)=0 & x(0)=0 & z(0)=0 & y(0)=0\\
p(1)=1.5 & x(1)=1 & z(1)=1 & y(1)=0\\
p(2)=1.5 & x(2)=1 & z(2)=0 & y(2)=0\\
p(3)=0 & x(3)=0 & z(3)=0 & y(3)=1\\
p(4)=0 & x(4)=0 & z(4)=0 & y(4)=0\\
\end{array} \right\},
\end{equation*}

\begin{equation*}
B=\left\{
\begin{array}{lccc}
p(0)=0 & x(0)=0 & z(0)=0 & y(0)=0\\
p(1)=0 & x(1)=0 & z(1)=0 & y(1)=0\\
p(2)=1.5 & x(2)=1 & z(2)=1 & y(2)=0\\
p(3)=1.5 & x(3)=1 & z(3)=0 & y(3)=0\\
p(4)=0 & x(4)=0 & z(4)=0 & y(4)=1\\
\end{array} \right\}.
\end{equation*}

On the other hand, let us observe that LP relaxation of Centralized-IP has another feasible point $C$, where
\begin{equation*}
C=\left\{
\begin{array}{lccc}
p(0)=0 & x(0)=0 & z(0)=0 & y(0)=0\\
p(1)=1 & x(1)=2/3 & z(1)=2/3 & y(1)=0\\
p(2)=1 & x(2)=2/3 & z(2)=0 & y(2)=0\\
p(3)=1 & x(3)=2/3 & z(3)=0 & y(3)=0\\
p(4)=0 & x(4)=0 & z(4)=0 & y(4)=2/3\\
\end{array} \right\}.
\end{equation*}
However, it is not possible to express point $C$ as convex combination of other two points. In other words, there is no $\alpha \in [0,1]$ satisfying $\alpha A + (1-\alpha)B=C$. Therefore, we conclude that the feasible region expressed by conv$(\X_i)$ does not include point $C$, thus it is tighter than the LP relaxation of the Centralized-IP.

\bibliographystyle{IEEEtran}

\bibliography{refs}

\end{document}